\documentclass[12pt]{article}
\usepackage{enumerate}
\usepackage{epsf,epsfig,amsfonts,a4wide}
\usepackage{amsfonts}
\usepackage{amsmath,amssymb,amsthm}
\usepackage{url}
\usepackage{amsmath,amssymb,amsthm}
\usepackage{amssymb}
\usepackage{amsthm}
\usepackage{epsf,epsfig,amsfonts,graphicx}
\usepackage{a4wide}
\newtheorem{Th}{Theorem}
\newtheorem{Lemma}{Lemma}

\newtheorem{proposition}{Proposition}

\newtheorem{Def}{Definition}

\theoremstyle{remark}
\newtheorem{Ex}{Example}
\newtheorem{Rem}{Remark}

\newcommand{\be}{\begin{equation}}
\newcommand{\ee}{\end{equation}}

\newcommand{\tr}{\mathrm{tr}\,}

\newcommand{\R}{\mathbb{R}}\newcommand{\Id}{\textrm{\rm Id}}

\newcommand{\trace}{\mathrm{tr}\,}\newcommand{\const}{\mathrm{const}}

\newcommand{\weg}[1]{}
\title{Splitting and gluing constructions 
 for geodesically equivalent pseudo-Riemannian metrics}
\author{Alexey V. Bolsinov\footnote{ School of Mathematics, 
 Loughborough University, 
 LE11 3TU, UK  \ \
 \quad {\tt A.Bolsinov@lboro.ac.uk} }
\and Vladimir S. Matveev\footnote{
Institute of Mathematics,
07737 Jena Germany  \ \ \quad {\tt  matveev@minet.uni-jena.de}} }
\date{} 
\begin{document}
\maketitle

\begin{abstract} Two metrics $g $ and $\bar g$ are geodesically equivalent, if they share  the same (unparameterized) geodesics.  We introduce two constructions that  allow one  to reduce many natural problems related to geodesically equivalent metrics, 
such as the classification of local normal forms  and the Lie problem (the description of projective vector fields),   
to the case when the $(1,1)-$tensor $G^i_j:= g^{ik}\bar g_{kj}$
 has one real eigenvalue, or two complex conjugate eigenvalues, and give first applications. As a part of the proof of the main result, we generalize Topalov-Sinjukov (hierarchy) Theorem for pseudo-Riemannian metrics

\end{abstract}

\section{Introduction}

\begin{Def} \label{geo}{\rm 
Let $g$ and $\bar g$ be Riemannian or 
pseudo-Riemannian metrics on the same   manifold $M^n$. We say that they are \emph{geodesically equivalent} (notation: $g\sim  \bar g$), if they have the same geodesics considered as unparametized curves.} \end{Def}

Given two metrics, we consider the $(1,1)-$tensor $L=L(g,\bar g)$ defined by
\begin{equation} \label{L}
L_j^i := \left(\frac{\det(\bar g)}{\det(g)}\right)^{\frac{1}{n+1}} \bar g^{ik}
 g_{kj}.\end{equation}
(Replacing  $g$  by $-g$ if necessary,  we can always assume  that \eqref{L} is well-defined.)

The goal of this paper is to give two constructions: gluing  and splitting. 

\begin{itemize}\item The  simplified  version of the  {\it gluing construction} does the following. 
Consider two manifolds $M_1$ and $M_2$ with  pairs of geodesically equivalent metrics $h_1\sim \bar h_1$ on $M_1$ and $h_2\sim \bar h_2$ on $M_2$. Assume  that the corresponding  $(1,1)$-tensor fields $L_1=L(h_1, \bar h_1)$ and $L_2=L(h_2, \bar h_2)$
 have no common eigenvalues in the sense that
for  any two points $x_1\in M_1$, $x_2\in M_2$ we have 
$$
\textrm{Spectrum}\, L_1(x_1) \cap \textrm{Spectrum}\, L_2(x_2) =\varnothing.
$$

Then one can naturally construct a pair of geodesically equivalent metrics $g\sim \bar g$ on the direct product   $M=M_1 \times  M_2$.  
These new metrics   $g$ and $\bar g$ differ from the direct product metrics  $h_1 + h_2$ and   $\bar h_1 + \bar h_2$ on $M_1\times M_2$, but can be obtained from  them by explicit  formulas involving $L_1$ and $L_2$
(see  (\ref {hh1}), (\ref{bh1}) below).  
We denote this by 
$$
(M,g,\bar g) = (M_1, h_1, \bar h_1) \times_{(L_1,L_2) } (M_2, h_2, \bar h_2)
$$
The corresposponding  $(1,1)-$tensor 
$L=L(g,\bar g)$ is, however,  the direct sum of $L_1$ and $L_2$  in the natural sense: 
for every 
$$
{\xi}= (\underbrace{\xi_1}_{\in T_{x_1}M_1}, \underbrace{\xi_2}_{\in T_{x_2}M_2}){\in T_{(x_1,x_2)}(M_1\times  M_2)} \ \textrm{ \ we have\ } \ L(\xi)= \left(L_1(\xi_1), L_2(\xi_2)\right).$$


\item The {\it splitting construction} is the inverse operation.  Its local version can be described as follows.  

Consider a manifold $M$ with two geodesically equivalent metrics  $g\sim \bar g$.
 Assume that at a point $p\in M$ 
  the corresponding $(1,1)-$tensor 
$L=L(g,\bar g)$  has at least two distinct eigenvalues that are not conjugates of each other. Then choose a partition of the spectrum of  $L$ into  two disjoint nonempty subsets 
$S_1\sqcup S_2$  such that each pair of complex-conjugate eigenvalues lies in the same subset.  Equivalently, one can say that  the characteristic polynomial  $\chi (t)$ of $L$  is factorised into two real polynomials $\chi_1(t) \cdot \chi_2(t)$ with no common roots  and  $S_i$ is just the set of roots of $\chi_i(t)$. 

Then there is  a neighborhood  $U=U(p) \subset M$ such that  in the above notation  the triple  
$(U, g, \bar g)$ can be presented as  
$$
(U,g,\bar g) = (U_1, h_1, \bar h_1) \times_{(L_1,L_2) } (U_2, h_2, \bar h_2),
$$
with $\dim U_1=\mathrm{deg}\, \chi_1$, $\dim U_2=\mathrm{deg}\, \chi_2$.  The splitting construction gives an explicit formula for the  pair of geodesically equivalent metrics $h_i \sim  \bar h_i$ on $U_i$  in terms  of $g, \bar g, L$ and the chosen  factorisation $\chi(t) =\chi_1(t) \cdot \chi_2(t)$  (see (\ref{h}), (\ref{bh}) below).   

\end{itemize}

The splitting construction is, of course,  more important:  it allows one to ``decompose''  $(M, g, \bar g)$ at least locally into simpler pieces $(U_i, h_i, \bar h_i)$ of the same kind. From the topological viewpoint  this decomposition (i.e.,  the existence of two transversal foliations of complementary dimensions on $M$)  is, in fact,  very natural and is  induced by the $(1,1)$-tensor field  $L$ itself.
We show that the Nijenhuis torsion of $L$ vanishes  (Theorem \ref{thm1}) and,  as result, the partition $S_1\sqcup S_2$  of the spectrum of $L$  leads immediately to two natural integrable distributions that define  on $M$ the desired  foliations (Theorem \ref{thm1b}).    The geometric part of the construction is much deeper.  We show in particular that  $M$ carries two hidden locally product  pseudo-Riemannian metrics $h$ and $\bar h$ which  can be  canonically reconstructed from $g$, $\bar g$  by means of  some non-trivial (and, in our opinion, beautiful) formulas (Theorem \ref{thm2}) involving the information about the partition of the spectrum.  Notice  that we formulate all of our results in invariant terms which makes it possible to apply them to  study global properties of manifolds with geodesically equivalent metrics.

In the Riemannian case, these constructions were obtained in \cite{hyperbolic,archive}, and played an important role  in recent developments   in the theory of geodesically equivalent Riemannian metrics.  Roughly speaking, these constructions allow  one to reduce many natural questions about geodesically equivalent metrics to the case when the $(1,1)-$ tensor 
 $L(g,\bar g)$ has one real eigenvalue, or two complex-conjugate eigenvalues.    

 We expect similar applications in the theory of  pseudo-Riemannian metrics.    The main idea  is still the same:   by iterating the splitting construction,  each pair of geodesically equivalent metrics  $(g ,  \bar g)$  can be canonically ``decomposed'' into  natural ``components''  $(h_i, \bar h_i)$,  $i=1,\dots, s$,  of smaller dimensions  with the same property  $h_i\sim \bar h_i$,   but with a much simpler structure of the tensor $L_i=L(h_i,\bar h_i)$.  Moreover,  some other important geometrical objects  appearing in this context inherit  this decomposition.  Thus, we are able to work with individual ``components'' and,  as a result, to reduce  essentially technical difficulties.  We give more details and  comments  in Sections \ref{beltrami}, \ref{lichnerowicz}.  
 
 It is worth noticing that compared to Riemannian geometry,  the gluing/splitting construction in the pseudo-Riemannian case seems to be even more important  because it replaces to some extent  the classical Levi-Civita theorem providing a canonical form for a pair of geodesically equivalent Riemannian metrics.  The absence of a reasonable analogue of the Levi-Civita theorem in the pseudo-Riemannian case is, in fact,  the main reason why
 the ``Riemannian'' proof  of the gluing/splitting theorems cannot be generalised to the case of  pseudo-Riemannian metrics, where one has to develop essentially different techniques (see 
 Remark \ref{riem}).

 A part of our  proof  of  Theorems 3 and 4 could be consider as  a separate  result. It is  a generalisation of a construction by Sinjukov \cite{sinjukov} and Topalov \cite{topalov} which makes it possible to build, starting  from a given pair of geodesically equivalent metrics $g\sim \bar g$,  a whole family  of geodesically equivalent pairs $g_f \sim \bar g_f$.  We give an invariant  explanation of this construction which works both in the Riemannian and pseudo-Riemannian context  and show that such a family is, in fact,  very large:  as  its parameter, $f$,  one can take an arbitrary real-analytic function satisfying certain natural assumptions,  see Section \ref{hijerarchy} for more details.

\subsection{Splitting construction}

\begin{Def} \label{lp}{\rm
A {\it local-product structure} on $M^n$ is a triple $(h, B_1, B_{2})$,
where  $h$ is a pseudo-Riemannian metric,  and $B_1$,
$B_{2}$ are  foliations of dimensions   $r$ and $n-r$ ($1\le r< n$)
such that in a neighborhood of each point $p\in M^n$ one can choose
local coordinates
$$
(x, y)= \bigr((x^1,x^2,...,x^r),(y^{r+1},y^{r+2},...,y^n)\bigl), 
$$
satisfying the following conditions:

1) the leaves of $B_1$ are given by $y=\const\in \R^{n-r}$; 

2) the leaves of $B_2$ are given by $x=\const\in\R^r$; 

3) the metric  $h$ takes the form
$$
ds_h^2 = \sum_{i,j=1}^r h_{ij}(x) \, dx^i dx^j + \sum_{\alpha,\beta=r+1}^n h_{\alpha\beta}(y)\,  dy^\alpha dy^\beta;
$$
in other words, the matrix of $h$ is block-diagonal  and its first
$r\times r$ block depends on the $x$-coordinates and
the second  $(n-r)\times (n-r)$ block depends on the $y$-coordinates:  
$$
h= 
\begin{pmatrix} h_1(x) & 0 \\  0 & h_2(y)\end{pmatrix}.
$$
}\end{Def}

\begin{Ex} 
A model example of manifolds with local-product structure is obviously
the direct product of two pseudo-Riemannian manifolds $(M_1^r, h_1)$ and
$(M_2^{n-r}, h_2)$.  In this case, the leaves  of the  foliations
$B_1$ and $B_2$ are respectivly $M_1^r \times \{y\}$, $y\in M_2^{n-r}$  
and $M_2^{n-r} \times\{ x\}$, $x\in M_1^r$. The  metric $h$ is the usual product metric $h_1 + h_2$.  
Locally, every  local-product structure is as in this model example.   
\end{Ex}

Geometrically, the existence of a local-product structure for a given metric $h$ is equivalent to the following condition:
$T_{p}M$ splits into the direct sum of two nontrivial orthogonal subspaces $U,V {\subset}T_{p}M$ invariant with respect to the holonomy group. Indeed, if $(h, B_1, B_2)$ is a local-product structure, then $U= T_{p}B_1$ and $V= T_{p}B_2$ are such subspaces. If the orthogonal  subspaces $U,V {\subset}T_{p}M$ (such that  $U\oplus V= T_{p}M$) are invariant with respect to the holonomy group, then, as  it was shown by De Rham  \cite{DR} and Wu \cite{Wu}, the parallel translation of these subspaces does not depend on the curve, and generates integrable distributions, whose integral manifolds form the foliations $B_1$ and $B_2$ satisfying  Definition~\ref{lp}.

For two metrics $g$ and $\bar g$,
   we consider the  $(1,1)-$tensor  $L=L(g,\bar g)$ given by \eqref{L}. 
Its characteristic polynomial will be denoted by
\begin{equation}
\chi:M\times \mathbb{R}\to \mathbb{R}, \  \  \chi(t):= \det(t\cdot \Id-L  ).
\end{equation} 
It is a  polynomial of degree $n$ whose coefficients are smooth functions on the manifold. Note that by our definition the polynomial is monic, i.e., its leading coefficient is $1$.

We say that a factorisation of $\chi(t)$ into two monic real polynomials 
$$
\chi (t) = \chi_1(t) \cdot \chi_2 (t), \quad   \mathrm{deg}\,\chi_i\ge 1, \quad \chi_i:M\times \mathbb{R}\to \mathbb{R} $$  
is {\it admissible} if $\chi_1(t)$ and $\chi_2(t)$ are coprime (i.e., have no common roots) at every point $x\in M$.

\begin{Ex} \label{ex1} 
Let $\lambda_1,\dots,\lambda_n$ be (possibly, complex) 
eingenvalues of $L$ at a point $p\in M$ 
counted with their  algebraic multiplicities.  Divide them into two nonempty groups (without loss of generality we think that $\lambda_1,\dots,\lambda_r$ lie in the first group and $\lambda_{r+1},\dots,\lambda_n$ in the second group) 
in such a way that $\lambda_i \ne \lambda_\alpha$  for  $i=1,\dots, r$, $\alpha=r+1,\dots,n$,  and   
pairs of complex-conjugate eigenvalues lie in the same group.
Then, the polynomials $$\chi_1(t):=(t-\lambda_1)\dots
(t-\lambda_{r}),\quad   \chi_2(t):=(t-\lambda_{r+1})\dots
(t-\lambda_{n}){}$$ are real and give an admissible factorisation of $\chi(t)$ at $p\in M$. 
By using the implicit function theorem, it is easy to see that this factorisation can always be extended onto some neighborhood $U(p)\in M$. 
\end{Ex} 

Locally every admissible factorisation is as in this example.  
However, the existence of a local admissible factorisation at every point does not imply its global existence.

Given an admissible factorisation
$\chi = \chi_1  \chi_2$, we consider the following two distributions $D_1, D_2$ on $M$:
\begin{equation}\label{D_i}  D_i=\mathrm{ker}\, \chi_i(L)\ , \ \ i=1,2 , 
\end{equation}
where $\textrm{ker}$ denotes the kernel. Here (and in all other places of the paper
where we consider polynomials in $L$ with coefficients being smooth functions on $M$) we treat the $(1,1)-$tensor $L$ as a linear operator acting on each tangent space, and a polynomial $f(L)$ in $L$  is the $(1,1)$-tensor of the form $f(L)=a_0(x) \cdot \mathrm{Id} + a_1(x) L + a_2(x) L^2 + \cdots + a_m(x) L^m$.  

It is easy to see that $D_1$ and $D_2$ are transversal 
distributions of complementary dimensions.    
Moreover, the distributions are invariant with respect to  $L$, and are 
mutually orthogonal with respect to the both metrics $g$ and $\bar g$.

\begin{Th} \label{thm1}
Let $g$ and $\bar g$ be geodesically equivalent pseudo-Riemannian metrics on $M$,  $L=L(g, \bar g)$ be the $(1,1)$-tensor associated with them,  
and  $\chi =\chi_1 \cdot \chi_2$ be an admissible factorisation of 
its characteristic polynomial.
Then, the distributions \eqref{D_i} are integrable.   
\end{Th}

By Theorem \ref{thm1}, the admissible factorisation implies the (local) existence of 
a coordinate system $(x^1,...,x^r,y^{r+1},...,y^{n})$ such that the distribution $D_1$ is generated by the coordinate vector fields 
$\partial_{x^i}$, $i=1,\dots,r$, and similarly $D_2$ is generated by  $\partial_{y^j}$, $j=r+1,\cdots, n$. Since the distributions are  invariant with respect to $L$, in this coordinate system the matrix of $L$ is block-diagonal:   
$$
L(x,y)=\begin{pmatrix} L_1(x,y) & 0 \\  0 & L_2(x,y)\end{pmatrix},
$$
where $L_1$ and $L_2$ are  $r\times r-$  and $(n-r)\times (n-r)-$ matrices respectively. 
 
 \begin{Th} \label{thm1b} Under the assumptions of Theorem {\rm \ref{thm1}} and in  notation above, the entries of $L_1$ depend on the $x-$variables only, and the   entries of  $L_2$ depend on the $y-$variables only,  so \begin{equation}\label{matl} 
L(x,y)=\begin{pmatrix} L_1(x) & 0 \\  0 & L_2(y)\end{pmatrix}.\end{equation}
 \end{Th}

By Theorem \ref{thm1}, the distributions $D_1 $ and $D_2$  generate two foliations  on the manifold. We denote them by $B_1$ and $B_2$, respectively.

\begin{Th}[Splitting construction] \label{thm2}
Let $g$ and $\bar g$ be geodesically equivalent pseudo-Riemannian metrics on $M$,  and
$B_1$, $B_2$ be the foliations generated by the distributions \eqref{D_i} related to an  
admissible factorisation $\chi =\chi_1 \cdot \chi_2$ of the characteristic polynomial of $L=L(g, \bar g)$.
Then,   the following two   $(0,2)-$tensors
 \begin{eqnarray}\label{h} h_{ij} &  = &  g_{ik}\left((\chi_2(L) +\chi_1( L))^{-1}\right)_{j}^k\\
 \bar  h_{ij} &  = &  \bar g_{ik}\left(\left(\frac{1}{\chi_2(0)}\chi_2(L) +\frac{1}{\chi_1(0)}\chi_1( L)\right)^{-1}\right)_{j}^k \label{bh}
\end{eqnarray}
are pseudo-Riemannian metrics  (i.e., symmetric and nondegenerate), and the triples $(h, B_1, B_2)$ and  $(\bar h, B_1, B_2)$  
are   local-product structures on $M$.    
Moreover, for every leaf $M_i$ of the foliation $B_i$,  the restrictions $h_i=h|_{M_i}$ and $\bar h_i|_{M_i}$  of the metrics 
$h$ and $\bar h$ to this leaf are geodesically equivalent,  and  
the tensor $L_i=L(h_i,\bar h_i)$ defined on $M_i$ by \eqref{L} coincides with the restriction $L|_{M_i}$ of $L$ to $M_i$, $i=1,2$.   Moreover, $\chi_1= \chi_{L_1}$ and $\chi_2= \chi_{L_2}$.

\end{Th}

In formula \eqref{bh},  the expression  $\chi_i(0)$,  $i=1,2$, denotes the smooth function on $M$ obtained by substituting $t=0$ into $\chi_i(t)$, namely,  $\chi_1(0) = (-1)^{r}\det L_1$, and $\chi_2(0)=(-1)^{n-r}\det  L_2$.

\begin{Rem} 
It might be  convenient to understand  the formulas (\ref{h}, \ref{bh})  in matrix notation: let us consider   the coordinate system
$(x^1,...,x^r,y^{r+1},...,y^{n})$ such that  the distribution $D_1$ is generated by the vectors $\partial_{x^i}$, and  the distribution $D_2$ is generated by the vectors $\partial_{y^j}$. 
Then, in this coordianate system  all matrices we use  to construct $h$ and $\bar h$ are 
block-diagonal with one $r\times r -$
and one $(n-r)\times (n-r)-$block:

 $$
g =\begin{pmatrix} g_1   & 0 \\  0 &   g_2 \end{pmatrix}\ , \ \  \bar g =\begin{pmatrix} \bar g_1   & 0 \\  0 &   \bar g_2 \end{pmatrix}\ , \ \ L =\begin{pmatrix} L_1   & 0 \\  0 &   L_2 \end{pmatrix}.$$
Hence, the matrices of $h$ and $\bar h$ must be block-diagonal as well, and direct calculations give us the formula

 \begin{equation}\label{math} 
h =\begin{pmatrix}   g_1   \chi_2(L_1)^{-1} & 0 \\  0 &  g_2   \chi_1(L_2)^{-1}\end{pmatrix}\ , \ \ \bar h =\begin{pmatrix} {\chi_2(0)}   \bar g_1 \chi_2(L_1)^{-1} & 0 \\  0 & {\chi_1(0)}  \bar g_2   \chi_1(L_2)^{-1}\end{pmatrix}
\end{equation}

\end{Rem}
\begin{Rem}
Although the geodesic equivalence relation is obviously symmetric, the tensor $L$  is not invariant with respect
to  the permutation of $g$ and $\bar g$. More precisely, the tensor
(\ref{L}) constructed for  $\bar g$, $g$ is the inverse of
(\ref{L}) constructed for  $g$, $\bar g$.  Moreover, two local product structures constructed from the pair of geodesically equivalent
 metrics $g,\bar g$ does not coincide with the local product structure constructed from the pair $\bar g$, $g$, though of course they are closely related, in particular  (after the proper choice of admissible factorisations) all these four local-product structures share the same foliations $B_1$ and $B_2$.   \end{Rem}

\subsection{Gluing construction}

Let $B_1$ and $B_2$ be two transversal foliations of complementary dimensions,   and   $h, \bar h$ be two metrics such that the triples $(h, B_1, B_2)$ and $(\bar h, B_1, B_2) $ are local-product structures.

The tangent spaces to the leaves of the foliations will be denoted  by $TB_1, TB_2$, and for every $p\in M$ we have $T_pM= T_pB_1 \oplus T_pB_2$. 

We denote by  $h_i, \bar h_i$   the  restrictions of the
 metrics $h, \bar h$  to  the leaves of  $B_i$,  $i=1,2$.  
The corresponding tensor \eqref{L} will be denoted by $L_i:= L(h_i, \bar h_i)$, and its characteristic polynomial by $\chi_i$, $i=1, 2$.
Assume in addition that the polynomials $\chi_1$ and   $\chi_2$ are coprime for all points of $M$.

Now,  consider the following two symmetric bilinear forms $g, \bar g$   on $M$:  for two tangent vectors
 $$u= (\underbrace{u_1}_{\in TB_1}, \underbrace{u_2}_{\in TB_2})\, , \  \ v=( \underbrace{v_1}_{\in TB_1}, \underbrace{v_2}_{\in TB_2}) \in TM $$
 we put  
 \begin{eqnarray} g(u,v) &  = &  h_1\left( \chi_2(L_1)( u_1), v_1\right)   + h_2\left(\chi_1(L_2)(u_2), v_2\right)  \label{hh1}  \\
  \bar g(u,v) &  = &  \frac{1}{\chi_2(0)}\bar h_1\left( \chi_2(L_1) (u_1), v_1\right)   + \frac{1}{\chi_1(0)}\bar h_2\left(\chi_1(L_2)(u_2), v_2\right).  \label{bh1}
\end{eqnarray}

\begin{Th}[Gluing construction] \label{thm3}

If $h_1$  is geodesically equivalent to $ \bar h_1$, and  $h_2$ is geodesically equivalent to $ \bar h_2$, then  the metrics $g,\bar g$ given by {\rm (\ref{hh1}, \ref{bh1})}  are geodesically equivalent too.
\end{Th}

\begin{Rem} 
It is an easy linear algebra to see that the tensor \eqref{L} constructed for the metrics (\ref{hh1}, \ref{bh1}) is the direct sum  of the tensors $L_1$ and $L_2$:  for every $u=(u_1, u_2)\in T_pB_1\oplus  T_pB_2=T_pM $ we have 
\begin{equation}
L(u)= (L_1(u_1), L_2(u_2)).
\end{equation}
\end{Rem}

\begin{Rem} 
It might again  be  convenient to understand  the formulas (\ref{hh1}, \ref{bh1})  in matrix notation: we  consider   the coordinate system
$(x^1,...,x^r,y^{r+1},...,y^{n})$ such that  $y-$coordinates are constant on the  leaves of  $B_1$ and  $x-$coordinates are constant on the  leaves of $B_2$. 
Then, in this coordianate system, the matrices  of $g$ and $\bar g$ are given by

 \begin{equation}\label{matg} 
g =\begin{pmatrix}   h_1   \chi_2(L_1) & 0 \\  0 &  h_2   \chi_1(L_2)\end{pmatrix}\ , \ \ \bar g =\begin{pmatrix} \frac{1}{\chi_2(0)}   \bar h_1 \chi_2(L_1) & 0 \\  0 & \frac{1}{\chi_1(0)}  \bar h_2   \chi_1(L_2)\end{pmatrix}
\end{equation}

\end{Rem}

\begin{Rem}

Comparing formulas \eqref{math} and   \eqref{matg}, we see 
that  the gluing construction  is inverse to the splitting. 
\end{Rem}

\begin{Rem}  \label{riem} If the case of Riemannian metrics,   Theorems \ref{thm2},\ref{thm3} were proven in \cite{hyperbolic,archive}. The proof is based on the Levi-Civita description of geodesically equivalent  metrics, which is not complete in the pseudo-Riemannian case, see the discussion in Section \ref{beltrami}.
 Below, we will show that the Levi-Civita description is an easy corollary of Theorems \ref{thm2},\ref{thm3}, see example  in Section \ref{beltrami}.\end{Rem}

 \subsection{Functions of  $(1,1)-$tensors and Topalov-Sinjukov \\ Theorem  for pseudo-Riemannian metrics } \label{hijerarchy}

Suppose a compact set   $K \subseteq \mathbb{C} $,  a function $f: {K} \to \mathbb{C} $, and a $(1,1)-$tensor $L$  (on $M$)  satisfy the following assumptions:

\begin{enumerate}[(i)] \item $\mathbb{C}\setminus K$ is connected (we do not require that $K$  is connected). \label{I}

 \item  $f:K\to \mathbb{C}$ is  a continuous function, and  the restriction $f |_{\mathrm{Int}\,K}$ is holomorphic. \label{ii}
\item $K$ is symmetric with respect to the $x-$axes: for every  $z\in K$ its conjugate $\bar z $ also lies in $K$. \label{iii}
 \item for every $z\in K$,  $f(z)= \bar f(\bar z)$, where the bar ``$\bar{ \hspace{1ex} }$"    denotes the complex conjugation.  \label{iiii}
  \item  $\mathrm{Spectrum}\, L(x)\subset \mathrm{Int}\,K$  for  every $x\in M$. \label{iiiii} \end{enumerate}

   Under the above assumptions (\ref{I}--\ref{iiiii}), one can naturally define  a 
$(1,1)-$tensor  $f(L)$: by  Mergelyan's theorem \cite{mergelyan}, the function $f$  can be uniformly  approximated by real polynomials $p_i$. 
We define $f(L)= \lim_{i\to \infty }p_i(L)$. It is an easy exercise  (see 
  for example   \cite[\S 1.2.2 -- 1.2.4]{higham}) to show that  the limit  exists,  is independent on the choice of the sequence $p_i$,   smoothly depends on $x\in M^m$ (actually, the function $f(L)$ is analytic in the entries of $L$), and behaves as a $(1,1)-$tensor when we change the variables.    
  
\begin{Ex} Polynomials $p(z)$ with real coefficients,    the functions  $e^z$, $\cos(z)$, $\sin(z)$ satisfy the above assumptions for every $L$, and we can naturally consider the standard operator functions $p(L)$, $e^L$, $\sin(L)$, $\cos(L)$ as smooth  $(1,1)$-tensor fields on $M$.
\end{Ex}   

\begin{Ex}\label{project} Let  $\chi =\chi_1 \cdot \chi_2$ be an admissible factorization of the characteristic polynomial of $L$  in a sufficiently small  neighborhood 
$ U(x_0)\subset M$  and   $D_1$, $D_2$ be the corresponding  distributions.  Then, the   natural projectors 
$P_i: TU \to D_i$ are functions of $L$ in the above sense.  Indeed, 
consider the corresponding partition  $S_1\sqcup S_2$ of  $\mathrm{Spectrum}\, L(x_0)$ into two disjoint subsets and let  $K_1$ and $K_2$  be $\varepsilon$-neighborhoods  of $S_1$ and $S_2$ respectively  ($\varepsilon$ is small enough so that  $K_1$ and $K_2$ do not intersect). Then, the function $$f_1:K\to \mathbb{C}, \  \ f_1(z):= \left\{\begin{array}{cc} 1 & \textrm{ for $z\in K_1$ } \\ 0 & \textrm{ for $z\in K_2$. } \end{array} \right. $$
satisfies the above assumptions, and we can define the function $f_1(L)$ which obviously coincides with  the projector $P_1$.  Similarly,  one can define $f_2(L)=P_2$.
\end{Ex}

 Notice that    
 for  every $x\in M $ we have: $\mathrm{Spectrum} f(L(x))=f(\mathrm{Spectrum}\,L(x))$.    In particular, if for  every  $x\in M$ we have $0\not\in f(\textrm{Spectrum}\,L(x))$, then the tensor $f(L)$ is nondegenerate.  
 If in addition $L$ is $g-$self-adjoint, then ${f(L)}$ is also $g-$self-adjoint, so 
  $g_f:= \left( g_{i\alpha}f(L)^\alpha_j \right)$ is a metric.

\begin{Th} \label{sintop}  Let $g$ and $\bar g$ be geodesically equivalent metrics and $L$  be the $(1,1)-$tensor given by \eqref{L}. Suppose $K\subset \mathbb{C}$  is a  compact  set, and $f:K\to \mathbb{C}$  is  a function such that $K,$ $ f$,  and    $L$ satisfy the assumptions (\ref{I} -- \ref{iiiii}) above. Assume in addition that for  every  $x\in M$ we have $0\not\in f(\mathrm{Spectrum}\, L(x))$.

 Then, the metrics $g_f:= \left( g_{i\alpha}f(L)^\alpha_j \right)$ and $\bar g_f:= \left( \bar g_{i\alpha}f(L)^\alpha_j \right)$ are also geodesically equivalent. 
\end{Th}

\begin{Rem} 
Partial cases of this theorem were proved by Sinjukov \cite{sinjukov} and Topalov \cite{topalov}. More precisely, Sinjukov proved this theorem assuming that the function $f(z)= z$. Topalov proved the theorem assuming the metric $g$ is Riemannian (in this case the eigenvalues of $L$ are real, and is sufficient to require that the function  $f$ is real-analytic).  
\end{Rem}

\section{History, motivation and possible applications} 
\label{history} 

The theory of  geodesically equivalent metrics  has a long and  fascinating history.  First  non-trivial examples 
 were discovered  by   Lagrange \cite{lagrange}. Geodesically equivalent metrics were studied by Beltrami \cite{Beltrami}, Levi-Civita \cite{LC}, Painlev\'e \cite{painleve} and other classics.   One can find more historical details in the surveys \cite{aminova1,mikes} and in the introductions to the papers \cite{hyperbolic,topology,archive}. 

The success of general relativity made necessary  to study geodesically equivalent pseudo-Riemannian  metrics. The textbooks \cite{eisenhart1,Petrov,Petrov2,hallbook} on pseudo-Riemannian metrics   have chapters on geodesically equivalent metrics. In the popular paper \cite{Weyl}, Weyl  stated  few interesting open   problems on geodesic equivalence of pseudo-Riemannian metrics.    Recent references (on the connection between geodesically equivalent metrics and general  relativity) include    Hall  et al \cite{hall,hall3,hall4}, Hall \cite{hall1,hall2}, Kiosak et al \cite{einstein}, Gibbons et al \cite{gibbons}.

 In the recent time, a huge progress was made in the theory  of  geodesically equivalent Riemannian metrics. The splitting/gluing constructions played a crucial role in this progress (as we mentioned in Remark \ref{riem}, in the Riemannian case, Theorems \ref{thm2}, \ref{thm3} were known). 
 We expect similar applications in the pseudo-Riemannian situation as well. The list of problems where  the  splitting/gluing  constructions  were used in the Riemannian case, and 
 are expected to be used in the pseudo-Riemannian case is below;   we discuss it  in detail in Sections \ref{beltrami}, \ref{lichnerowicz}.

 {\bf Beltrami Problem\footnote{ Italian original from \cite{Beltrami}: \emph{
La seconda $\dots$  generalizzazione $\dots$ del nostro problema,   vale a dire:   riportare i punti di una superficie sopra un'altra superficie in modo  che alle linee geodetiche della prima corrispondano linee geodetiche della seconda}}.} \emph{ Describe all pairs of geodesically equivalent metrics.}

{\bf Lie  Problem\footnote{ German original from \cite{Lie}, Abschn. I, Nr. 4, \\
\emph{Man soll die Form des Bogenelementes einer jeden Fl\"ache bestimmen,
deren geod\"atische Kurven eine  infinitesimale Transformation gestatten}}.}
\emph{Find all metrics $g$ admitting infinitesimal projective transformations.}

{\bf Lichnerowicz Conjecture\footnote{the attribution to Lichnerowicz is folkloric; we did not find a paper of Lichnerowicz where he states  this conjecture. Certain papers refer to   this statement  as to  a classical conjecture,  see the discussion in  \cite{archive}.  Lichnerowicz conjecture answeres the question asked by Schouten in \cite{Schouten}}.} \emph{  Let a connected  Lie group $G$  act on a complete
 connected  manifold $(M^n, g)$ of dimension
 $n\ge 2$   by projective
transformations.  Then, it  acts by affine transformations, or  for some $c\in \mathbb{R}\setminus \{0\}$ the metric $c\cdot g$
is the Riemannian  metric of  constant positive sectional  curvature $+1$.}

Recall that a \emph{projective transformation} of a Riemannian manifold is a diffeomorphism of the manifold that takes unparameterized geodesics to geodesics. Local  projective transformations obviously form a Lie pseudo-group, its generators, i.e., vector fields whose local flow takes unparameterized geodesics to geodesics are called by Lie   \emph{infinitesimal projective transformations}. In the modern terminlogy,  they are called \emph{projective vector fields}; we will use this terminology in the paper.

 \subsection{Motivation I. Normal form of geodesically equivalent metrics:  Beltrami problem  } \label{beltrami}

If the eigenvalues of $L$ do not bifurcate at  a point (this condition is fulfilled almost everywhere), the answer to Beltrami's question was given by Levi-Civita \cite{LC} under the assumption that  the metrics are  Riemannian.

A local description of geodesically equivalent pseudo-Riemannain metrics, which might be treated as a pseudo-Riemannian analog of the Levi-Civita theorem,   is considered to be done by Aminova \cite{Aminova}.  Unfortunaltely, the authors of the present paper do not understand her result (we do not doubt that the result is  correct).   Moreover,  we have checked that 
 in all 6 papers which refer to \cite{Aminova} according to MathSciNet, the authors cited Aminova's theorem to give an  overview of the subject only, but did not really use it.

 The statement of the 
main theorem of \cite{Aminova} is on two pages, and one more page is devoted to explanation of the formulas in the theorem.    This is probably the reason why this result is hardly applicable and was, to the best of our knowledge, never used.

The splitting and gluing constructions suppose to make the description of geodesically equivalent pseudo-Riemannian metrics much simpler. Indeed, at almost every point  $p\in M$ the eigenvalues of $L$ do not bifurcate,  i.e.,  the algebraic multiplicity of each eigenvalue $\lambda_i$  is locally constant  and is equal to $k_i$.  In a small neighborhood of  such a point,  the eigenvalues (both real and complex) can be treated as  smooth functions of~$x$, and we can factorise the characteristic polynomial of $L$  into  a product 
$\chi = \chi_1 \cdot  \ldots  \cdot  \chi_m$ of polynomials  $\chi_i:\mathbb{R}\times M\to \mathbb{R}$ of two kinds:  either $\chi_i(t)= (t -\lambda_i)^{k_i}$ for a real eigenvalue $\lambda_i$,  
or $\chi_i(t)= (t -\lambda_i)^{k_i}(t -\bar \lambda_i)^{k_i} $ for a pair of complex conjugate eigenvalues $\lambda_i, \bar \lambda_i$.

Repeatedly applying the splitting construction $m-1$ times in a small neighborhood $U(p)$, we can construct   
 local coordinates $x_1, \dots , x_m$,  where $x_i=(x_i^1, \dots, x_i^{l_i}) \in \mathbb R^{l_i}$ 
 ($l_i = k_i$ if $\lambda_i$ is real, and $l_i = 2k_i$   
  for a couple of complex conjugate eigenvalues $\lambda_i, \bar\lambda_i$),   
and pairs of geodesically equivalent metrics  $h_i (x_i) \sim \bar h_i (x_i)$
such that

1)  each eigenvalue $\lambda_i$ of $L$ depends on $x_i$ only: 
$\lambda_i=\lambda_i(x_i)$,

2)  the characteristic polynomial of $(1,1)$-tensor   $L_i = L(h_i, \bar h_i)$ associated with $h_i$ and $\bar h_i$ by means of (\ref{L})   is exactly $\chi_i(t)$,

\weg{\item the $(1,1)$-tensor   $L_i = L(h_i, \bar h_i)$ associated with $h_i$ and $\bar h_i$ by means of (\ref{L})  
has one real eigenvalue $\lambda_i(\bar u_i)$ for $i=1,\dots, m$,  or a couple of complex conjugate eigenvalues $\lambda_i(u_i), \bar\lambda_i(u_i)$, for $i=m+1, \dots$
in particular,  the characteristic polynomial of  $L_i$ is exactly $\chi_i(t)$,}

3)  in this coordinate system, $g$ and $\bar g$ simultaneously take the following  block diagonal form 
\begin{equation}
g=\begin{pmatrix}
h_1 \hat \chi_1 (L_1)  & & \\
& \!\!\!\!\!\! \ddots & \\
& & \!\!\!\!\!\! h_m\hat \chi_m (L_m)
\end{pmatrix},  
\quad 
\bar g=\begin{pmatrix}
 \frac{1}{\hat \chi_1(0)}
\bar h_1 \hat \chi_1 (L_1)  & & \\
&\!\!\!\!\!\! \ddots & \\
& & \!\!\!\!\!\! \frac{1}{\hat \chi_m(0)}\bar h_m \hat\chi_m (L_m) 
\end{pmatrix}
\label{decomposition}
\end{equation}
where $\hat \chi_i (t) = \prod _{j\ne i} \chi_j (t)$.

Thus we have
\begin{proposition}
\label{decomp}
In a neighborhood of a regular point $p\in M$,  the geodesically equivalent metrics $g$ and $\bar g$  can be simultneously reduced by an appropriate choice of local coordinates $x_1, \dots , x_m$,  $x_i\in \mathbb R^{l_i}$,  to the form \eqref{decomposition}, where   $h_i$ and $\bar h_i$ are geodesically equivalent metrics depending on $x_i$-coordinates only  and such that the corresponding  $(1,1)$-tensor 
$L_i=L(h_i, \bar h_i)$ has either one single real eigenvalue $\lambda_i(x_i)$, or  a single pair of complex conjugate eigenvalues $\lambda_i(x_i), \bar\lambda_i(x_i)$.

\end{proposition}

Thus, in order to describe geodesically equivalent metrics near a generic point it is sufficient to do this under the assumption that $L$ has one real eigenvalue, or two complex-conjugate eigenvalues.  In the Riemannian  case,   we will illustrate this idea  by  the following  

{\bf Example:  Levi-Civita Theorem follows from  the splitting construction.}
Let the geodesically equivalent metrics $g\sim \bar g$ be  Riemannian. 
Then all eigenvalues of $L$ are real and positive, and $L$ is semi-simple. 
Hence, the tensor $L(h_i, \bar h_i)$ for geodesically equivalent  metrics $h_i\sim \bar h_i $ discussed above 
 is $\lambda_i\cdot \mathrm{Id}$. If the multiplicity $k_i$ of $\lambda_i $ is $\ge 2$, then by the classical result of Weyl \cite{Weyl2} the metrics $h_i$ and $\bar h_i$ are proportional, i.e., 
 $\bar h_i:= \frac{1}{\lambda_i^{k_i+1}}h_i$, where the eigenvalue $\lambda_i$ is a constant.  
 
 If  the multiplicity $k_i$ of $\lambda_i $ is $1$, then $h_i$  is one-dimensional and we can obviously choose  $x_i$ in such a way that   $h_i= dx_i^2$ and $\bar h_i= \frac{1}{\lambda_i(x_i)^2}dx_i^2$.  
 
 Without loss of generality we can assume that the first $r$ eigenvalues $\lambda_1,...,\lambda_r$ have multiplicity 1, and the last $m-r$ eingenvalues have multiplicity $\ge 2$.  Then, in the chosen coordinate system,  the direct product metrics  $h=h_1+h_2+\cdots + h_m$ and $\bar h=\bar h_1+\bar h_2+\cdots + \bar h_m$ are given by 
 $$
\begin{array}{cccc} ds_h^2 &  = &  \sum_{i=1}^r{dx_i}^2 &+   \sum_{i=r+1}^m \left[\sum_{\alpha_i, \beta_i=1}^{k_i}(h_i(x_i))_{\alpha_i \beta_i}dx_i^{\alpha_i}dx_i^{\beta_i} \right]\\ 
ds_{\bar h}^2  &  = &  \sum_{i=1}^r\frac{1}{\lambda_i(x_i)^2}{dx_i}^2 &+   \sum_{i=r+1}^m  \left[\frac{1}{\lambda_i^{k+1}}\sum_{\alpha_i, \beta_i=1}^{k_i}(h_i(x_i))_{\alpha_i \beta_i} dx_i^{\alpha_i}dx_i^{\beta_i}\right]. 
\end{array}
 $$
 
 Here the functions $\lambda_i$ are constant for $i>r$ and depend only on the corresponding variable $x_i$ for $i\le r$. The metrics $ h_i$, $i=r+1,...,m$ can be arbitrary, but  their  entries  $(h_i)_{\alpha_i\beta_i}$  must depend on the coordinates $x_i= (x_i^1,...,x_i^{k_i})$ only.

 Applying Proposition \ref{decomp} and formula \eqref{decomposition}, we obtain for $g$ and $\bar g$ the following  form:
 
 $$
\begin{array}{cccc} ds_g^2 &  = &  \sum_{i=1}^rP_i{dx_i}^2  & + \sum_{i=r+1}^m \left[P_i \sum_{\alpha_i, \beta_i=1}^{k_i}(h_i(x_i))_{\alpha_i \beta_i}dx_i^{\alpha_i}dx_i^{\beta_i}\right] \\ 
ds_{\bar g}^2  &  = &   \sum_{i=1}^r P_i \rho_i {dx_i}^2 & + \sum_{i=r+1}^m \left[ P_i \rho_i\sum_{\alpha_i, \beta_i=1}^{k_i}(h_i(x_i))_{\alpha_i \beta_i} dx_i^{\alpha_i}dx_i^{\beta_i}\right], 
\end{array}$$
 where 
 \begin{equation}
 P_i:= \pm \prod_{j\ne i} (\lambda_i- \lambda_j), \  \  \ \rho_i:= \pm \frac{1}{\lambda_i \, \prod_{\alpha} \lambda_\alpha}. 
 \end{equation}
(the signs $\pm$ should be chosen so that all  $P_i$ and $\rho_i$ are  positive). 
This is precisely the Levi-Civita normal form from \cite{LC} for geodesically equivalent Riemannian 
metrics!

As we mentioned above, the results by  Levi-Civita and Aminova hold in a neighborhood of almost every point.  More precisely, 
such a  point  $p\in M$, which we call {\it regular}, is characterised by the property that the structure of the Jordan normal form of $L$
  (the number of  Jordan blocks and their dimensions) is the same for all points in some neighborhood $U(p)$.
  Regular points form an open everywhere dense subset of $M$. The other points will be called \emph{singular}.

For global questions (in particular, for the description of the topology of closed 
manifolds admitting geodesically equivalent metrics),  it is  also necessary to solve Beltrami problem near  singular points.   
     For the Riemannian case,  it was  done in  \cite{hyperbolic,archive,oshemkov}.  
 The result essentially used the  splitting construction: arguing as above, the problem was reduced to  a few cases of  simple bifurcations, which were considered separately.  We expect the same application of the splitting construction in the pseudo-Riemannian case.

 \subsection{Motivation II. Projective transformations: Lie problem  and  Lichnerowicz conjecture} 
 \label{lichnerowicz}
  For Riemannian manifolds,  splitting construction  found  an  important application in the theory of projective transformations.  Projective transformations is a very classical object of study. The first examples are due to Beltrami \cite{Beltrami}; as we mentioned at the beginning of Section \ref{history},  the problem of local description of projective transformations was explicitely stated by Lie \cite{Lie}. In the global setting, namely under the assumption that  $(M,g)$ is complete, the problem   was formulated by Schouten \cite{Schouten};  the Lichnerowicz conjecture mentioned above is a hypothetical solution of this problem. 
   
   In this section we  assume that $\dim M\ge 3$. The reason for this is that 
    in dimension 2 the Lie problem was solved in \cite{bryant,alone}, where a complete list of local metrics admitting projective vector fields was constructed. The list is pretty simple (explicit formulas involving  only  elementary functions) and it should  not be  very complicated to understand which metrics from this list can be prolonged  up to a complete  metric.

  As it was observed by Fubini and Solodovnikov,  
 in dimension  $\ge 3$, the degree of mobility plays a crutial role in the description of 
     projective transformations of a given manifold. Recall that \emph{the degree of mobility}
     of a  metric $g$ is the dimension of the space  of the solutions of the equation \eqref{main} considered as an equation on $L$. The degree of mobility has a clear geometrical meaning: locally, it coincides with the  dimension of the 
     set   of metrics geodesically equivalent to $g$ equipped with the natural topology. 
     
      In particular,  the condition ``degree of mobility = 1''   means that  $g$ does not admit any  geodesically equivalent ``partner'' except for  $\bar g=\const\cdot   g$. In this case projective transformations  are just homotheties of $g$.

Riemannian     metrics (on ${M}^{n\ge 3}$)   with  degree of mobility greater or equal to three, and their projective vector fields  were  locally described by Solodovnikov \cite{sol}; his result was improved  by Shandra \cite{shandra}, who  obtained  a much shorter  description based on a certain trick. The results of   
 \cite{kiosak2} show that 
 \begin{itemize} 
 \item the metrics with degree of mobility $\ge 3$ always admit projective vector fields, 
 \item  the trick  used in \cite{shandra} survives for  pseudo-Riemannian metrics as well, and  gives 
 a description of all metrics with degree of mobility  $\ge 3$, and their projective vector fields.   
     \end{itemize}

     Thus, in order to get the complete  solution of  the Lie Problem in dimension $\ge 3$, it is sufficient to consider metrics with degree of mobility $2$.  It turns out that under this assumption,  the problem can be reduced to the analysis of a  certain 1st order system of PDEs.  This system is universal in the sense that it does not depend on  the metric and remains the same if we pass from the whole manifold $M$ to the  components of smaller dimensions  obtained by means of the splitting construction.  We describe this  idea in brief below.

First of all we notice that in all our considerations a pair $g \sim \bar g$ of geodesically equivalent metrics can be replaced by the pair $(g, L)$ where $L$ is the (1,1)-tensor  defined by  \eqref{L}.  Indeed, $\bar g$ can be uniquely reconstructed from $g$ and $L$  as
$\bar g=    
\frac{1}{\det(L)} g L^{-1}$.  We shall say that $L$  and  $g$ are {\it compatible},  if  $g$ and 
$\bar g=    
\frac{1}{\det(L)} g L^{-1}$ are geodesically equivalent. 
The $(1,1)$-tensors $L$ compatible with $g$ form a finite-dimensional vector space whose dimension is exactly the degree of mobility of $g$.

Now let $v$ be  a projective vector field for the metric $g$  (we work in a small connected neighborhood $U$ of a point $p\in M$).  Consider the tensor 
\begin{equation}\label{apl}
\tilde L:= g^{-1} {\mathcal{L}}_v g - \tfrac{1}{n+1} \trace( g^{-1} {\mathcal{L}}_v g) \cdot \Id
\end{equation}
where   ${\mathcal L}_v$ denotes  the Lie derivative with respect  to $v$.

\begin{Lemma}[\cite{archive}]  \label{ap1}  
The tensor $\tilde L$ given by \eqref{apl}  is compatible with $g$. 
\end{Lemma}  

Lemma \ref{ap1} is an easy corollary of the compatibility condition \eqref{main} below, its 
 proof can be found for example in \cite[Section 2.1]{archive}, see Theorem 3 there. In an equivalent form the statement appeared already in Fubini \cite{Fubini1} (for dimension 3) and Solodovnikov \cite{sol} (under additional assumptions).

Now assume that the degree of mobility of $g$ is 2, so that the vector space of admissible $(1,1)$-tensors is two-dimensional.  Clearly, the identity tensor $\Id$ 
is admissible and we can take it as the first basis vector. We choose the second one, denoted by  $L$,  in such a way that  $\det L\ne 0$  in some neighborhood of $p\in  M$.

\begin{Lemma}\label{ap2}  Under the above assumptions,  the tensor $L$ and the projective vector field $v$ satisfy the following equation 
\begin{equation} \label{lil}
 {\mathcal L}_v L = \alpha \cdot L^2 + \beta\cdot L + \gamma \cdot \Id,  
\end{equation}
where  $\beta$ is a function given by {\rm  $ \beta := -{\mathcal L}_v (\log|\det L|)\cdot L +\alpha \cdot \trace(L)+ \gamma \cdot \trace(L^{-1}) + \nu$},  and $\alpha, \beta, \nu$  are some constants.
\end{Lemma} 

{\bf Proof.} Consider the metric $\bar g=\frac{1}{\det(L)} g L^{-1}$. 
This metric is geodesically equivalent to $g$ and therefore
$v$ is a projective vector field also for $\bar g$.  Clearly, 
the degree of mobility of $\bar g$ coincides with that of $g$. Hence, the space of  $\bar g$-compatible
$(1,1)$-tensors  is two-dimensional 
and is generated by $\mathrm{Id}$ and $L^{-1}$. 

Thus, applying Lemma \ref{ap1} to $g$ and $\bar g$,  we have:
$$
\begin{aligned}
g^{-1} {\mathcal{L}}_v g - \tfrac{1}{n+1} \trace( g^{-1} {\mathcal{L}}_v g)&=a\cdot L + b\cdot\mathrm{Id},  \\ 
\bar g^{-1} {\mathcal{L}}_v \bar g - \tfrac{1}{n+1} \trace( \bar g^{-1} {\mathcal{L}}_v \bar g)&= c\cdot L^{-1} + d\cdot\mathrm{Id}, 
\end{aligned}
$$
where $a,b,c,d\in\R$ are some constants. It is easy to see that these equations can be rewritten in the following equivalent form:
\begin{eqnarray} 
g^{-1} {\mathcal{L}}_v g &=&a\cdot L + \bigl(a\cdot \trace (L) + (n+1)b\bigr)\cdot\mathrm{Id},  \label{somen1}\\ 
\bar g^{-1} {\mathcal{L}}_v \bar g &=&c\cdot L^{-1} + \bigl(c\cdot \trace(L^{-1}) + (n+1)d\bigr) \cdot\mathrm{Id}, \label{somen2} 
\end{eqnarray}

On the other hand, we have
$$
\begin{aligned}
\bar g^{-1} {\mathcal{L}}_v \bar g &= \det L\cdot  L g^{-1} {\mathcal L}_v\left( \frac{1}{\det L} g L^{-1} \right)  = \\   
&=-{\mathcal L}_v(\log|\det L|)\cdot \Id + L \, g^{-1} {\mathcal L}_v g \, L^{-1} - {\mathcal L}_v (L)\cdot L^{-1}
\end{aligned}
$$
Substituting \eqref{somen2} to the left hand side and \eqref{somen1} to the right hand side (into the middle term) of this relation, we get  
$$
\begin{aligned}
-{\mathcal L}_v(\log(|\det(L)|))\cdot \Id   + 
a\cdot  L + (a \cdot \trace(L) + (n&+1)b)\cdot \Id  - {\mathcal L}_v (L)\cdot L^{-1} = \\ &= c L^{-1}  +   (c \cdot \trace(L^{-1})+ (n+1)d)\cdot \Id . 
\end{aligned}
$$

After multiplication by $L$, this relation can be rewritten as
$$   {\mathcal L}_v (L) = 
a\cdot  L^2  + \left(-{\mathcal L}_v(\log |\det L|) + a \, \trace L- c  \, \trace L^{-1} + (n+1)(b - d)\right) \cdot L    - c\cdot \Id,  
$$   which is equivalent  to  \eqref{lil}. \qed

\begin{Lemma}  \label{ap3}  Let $\chi= \chi_1 \cdot \chi_2$ be an admissible factorisation of the characteristic polynomial of $L$ and
$(x^1,...,x^r,y^{r+1},...,y^n)$ be the local coordinate system induced by this factorisation (see Theorems \ref{thm1}, \ref{thm1b}). 
Then the projective vector field $v$ splits in this coordinate system in the sense that 
$$
v= \sum_{i=1}^r v^i(x)\frac{\partial }{\partial x^i} + \sum_{j=r+1}^n v^j(y)\frac{\partial }{\partial y^j},
$$
i.e., the first $r$ entries of $v$  do not depend on the $y-$coordinates, 
and the last $n-r$ entries of $v$  do not depend on the $x-$coordinates. 
Moreover, the function $\beta$ from Lemma \ref{ap2} is constant.
\end{Lemma}

{\bf Proof.} By Theorems \ref{thm1}, \ref{thm1b},  in the coordinate system $(\bar x,\bar y)$  the matrix of $L$ has the block-diagonal form \eqref{matl}. 
Then, the matrices of $L^2$ and of  ${\mathcal L}_v(\log |\det L|) \cdot L$ are block-diagonal as well with the same dimensions of blocks implying that ${\mathcal L}_v(L)$ is also block-diagonal. Then, the first $r$ components of $v$ depend on $x-$coordinates
 only, and the last $(n-r)-$components of $v$ depend on $y-$coordinates only: $ v=(v_1(x), v_2(y))$ where $v_1\in TB_1 $
 and  $v_2\in TB_2$. Multiplying the equation \eqref{lil1} by $L^{-1}$ and rewriting it in the form 
 $$
 ({\mathcal L}_v L)\cdot L^{-1} -  \alpha \cdot L - \gamma \cdot L^{-1}  =  \beta \cdot \Id 
$$
 we see that the first block of the left-hand side is independent of $y$ and the second block of the left-hand side is  independent of $x$ implying $ \beta =\const$. \qed

Thus, under the  additional assumption that an admissible factorisation exists,  the tensor $L$ and the vector field $v$ satisfy the following equation 
\begin{equation} \label{lil1}
{\mathcal L}_v ( L)  = \alpha \cdot L^2 + \beta\cdot L + \gamma \cdot \Id,  
\end{equation}
where $\alpha$,  $\gamma$ and $\beta$  are certain constants.  Moreover, in the notation of Theorem \ref{thm1b}, the $r\times r-$ resp. $(n-r)\times (n-r)-$ matrices $L_1$ and  $L_2$  viewed as $(1,1)-$tensors  satisfy the  
equations  \begin{eqnarray} 
{\mathcal L}_{v_1} ( L_1)  &= &  \alpha \cdot L_1^2 + \beta\cdot L_1 + \gamma \cdot \Id,  \label{un1}\\ 
{\mathcal L}_{v_2} ( L_2)  & =& \alpha \cdot L_2^2 + \beta\cdot L_2 + \gamma \cdot \Id, \label{un2} \end{eqnarray}
where the components of the vector $v_1$ are the   first  $r$ entries of $v$, and  components of the vector $v_2$ are the last $n-r$   entries of $v$. 
In other words, equation \eqref{lil1}  splits into two  independent  equations (\ref{un1},\ref{un2}). 

We can go further: if $\chi$ is a product of several mutually prime monic polynomials  $\chi_1,\dots,\chi_m$, we can split the equation (\ref{un2}) into $m$ equations of the similar form by repeating the above construction. After finitely many steps, we land at an   independent system of PDE of the form 
$$
{\mathcal L}_{v_i} ( L_i)  =   \alpha \cdot L_i^2 + \beta\cdot L_i + \gamma \cdot \Id,   \  \  \ i=1,...,m
$$ 
where  each $L_i$ and $v_i$ depend on the corresponding coordinates only.  Moreover, in a neighborhood of a regular point 
we may assume that $L_i$ has  one real eigenvalue, or two complex-conjugate eigenvalues.

  Thus, it is sufficient to solve equation \eqref{lil1} under the assumption that $L$ has one real eigenvalue, or two complex-conjugate eigenvalues.
 
In the Riemannian case, equation \eqref{lil1}  was obtained by other methods  and  played an important role in the proof of the projective Lichnerowicz conjecture \cite{archive}, and in the local description of projective vector fields \cite{Fubini1,sol,bryant,alone}. This equation  is relatively simple and  can be  solve in the Riemannian case explicitly.  Analysis of its solutions was, in fact, one of the most principal steps in the solution of the Lie and Schouten  problems
(under the additional assumption that the degree of mobility is at most two). 
       
      We expect similar applications in the pseudo-Riemannian case. If $L $ is
      semi-simple and has real eigenvalues, then at least locally there is no essential difference  from the  Riemannian case.  The possible difficulties might appear if $L$  has   Jordan blocks or  complex-conjugated eigenvalues. They  appear  already in dimension 2 (though in dimension 2 the difficulties  have been overcome):  two dimensional metrics 
     admitting projective vector fields 
      in the case when  $L$ is a Jordan block are given by much more complicated formulas than those which appear for $L$  semi-simple, see \cite[Theorem~1]{alone}.

\subsection{Motivation  III.  Topological  geodesic rigidity problem} 
\label{rigidity}

A local version of the splitting/gluing constructions is sufficient for local problems such as  finding normal forms for a pair of geodesically equivalent metrics and the Lie problem. For global problems (when the underlying manifold is assumed to be closed or/and complete) 
such as the Lichnerowicz conjecture discussed above, a general version of the splitting/gluing constructions was useful in the Riemannian case, and is expected to be useful in the pseudo-Riemannian one. Another example of global problems is 

{\bf Topological geodesic rigidity problem.} {\it On what manifolds, the unparametrized geodesics of any metric $g$  detemine this metric uniquely (up to multiplication by a constant)?}  

In other words, we want to know whether or not a given manifold $M$ admits at least one pair $g \sim \bar g$ of non-proportional geodesically equivalent (pseudo)-Riemannian metrics. 

All results in this area are, in fact, very recent: the first obstructions that prevent a manifold to possess such metrics were found in \cite{MT} (it was proved that geodesically equivalent metrics on a closed surface of genus $\ge 2$ are always proportional).
Moreover, it was  generally believed and 
explicitly stated in the survey \cite{mikes}  that it is hard to obtain such obstructions. Main   research was concentrated in the  following direction: find out   what metrics  are \emph{geodesically rigid}, in the sense they 
{ do}  not possess  nontrivial  geodesically equivalent ``partners''. 
 One of the first results in this direction is \cite{Si2}:  every metric geodesically equivalent to an irreducible symmetric metric of nonconstant curvature  is proportional to it. 
 Local geodesic rigidity problem was very popular in 60th--80th; there are more than 100 papers devoted to it, see surveys \cite{mikes,aminova1}.

Later, one started to investigate  the 
global geodesic rigidity problem (assuming that $M$ is closed or $g$ is complete).  A typical result is as follows: if $g$ is  a Riemannian  Einstein irreducible metric  of nonconstant 
curvature  on a closed manifold, then it is geodesically rigid \cite{mikes,Mi2} (though locally there exist Einstein Riemannian metrics of nonconstant curvature  that are not geodesically rigid).  The most standard (de-facto, the only) way  to prove 
such results was to use tensor calculus to 
 canonically obtain a nonconstant  function $f$ 
  such that $\Delta_g f= \const\cdot  f$ with $\const \ge 0$, which of course cannot exist on a closed Riemannian manifold.

Here is one of examples that shows how the splitting/gluing construction works in global setting without any specific assumptions about the metric (we simply combine Theorem~\ref{thm2} with the results 
by Wu \cite{Wu} and De Rham \cite{DR}).

\begin{proposition}
\label{global}
Let $M$ admit a pair of geodesically equivalent (pseudo)-Riemannian metrics $g$ and $\bar g$,  one of which is complete. Assume that
the characteristic polynomial $\chi$ of $L=L(g,\bar g)$ possesses an admissible global factorisation $\chi = \chi_1 \cdot \chi_2$ on $M$. Then $M$ admits a (pseudo)-Riemannian metric $h$ with a reducible  holonomy group  and its universal cover $(\widetilde M, \pi^*h)$  (where $\pi : \widetilde M \to M$ is the natural projection)  is  the  direct product of 
two (pseudo)-Riemannian manifolds $(M_1, h_1)$ and $(M_2, h_2)$ such that the fundamental group $\pi_1(M)$ acts on $M_1\times M_2$ by fiberwise isometries.
\end{proposition}

 Obviously, the above property is a very strong topological restriction on $M$. For example, in dimension 2,
 there are only five connected  manifolds of this kind: $\R^2$,  torus,  Klein bottle, M\"obius strip and cylinder $\R^1\times S^1$.
 
  The additional assumption in Proposition~\ref{global} about 
 the existence of a global admissible factorisation is, of course, very essential. There are many examples of geodesically equivalent metrics where this condition fails.
It appears, however, that in the Riemannian case the nonexistence of such a factorization implies that the fundamental group of $M$ is finite 
(assuming that $M$ is closed and admits two non-proportional geodesically equivalent metrics) \cite{hyperbolic}.

Combining these two observations, it was possible to describe all geodesically rigid 3-manifolds \cite{topology}, and to prove that topologically-hyperbolic manifolds are geodesically rigid, in the sense that two geodesically equivalent metrics on such  manifolds are proportional \cite{hyperbolic}. \weg{ The new methods of \cite{hyperbolic}, \cite{topology}, which lead to these topological results in the Riemannian case, 
essentially used the gluing/splitting construction. }

At the present point, it is not clear whether the  second observation could be generalized to the pseudo-Riemannian case, i.e., whether the nonexistence of an admissible factorisation implies that the  fundamental group is finite. However, the known examples allow us to suggest
that compared to the classical Riemannian situation, the topological obstructions 
to the existence of geodesically equivalent metrics are perhaps even stronger in the pseudo-Riemannian case. We conclude this section by the following     

{\bf Conjecture.} {\it   Let $g\sim \bar g$ be geodesically equivalent non-proportional  metrics  on 
 a closed 3-dimensional manifold  $M^3$. If at least one of the metrics   
  is complete,  and at least one of the metrics   has signature $(-,+,+)$,  then  $M$ is a Seifert manifold with zero Euler class.  }

\section{Proofs} 
\subsection{Vanishing of Nijenhuis torsion and  proof of Theorems~\ref{thm1}, \ref{thm1b} }

Throughout the paper we shall use the following analytic condition for two metrics $g$ and $\bar g$ to be  geodesically equivalent.

\begin{proposition}\label{maineq} Let $g$, $\bar g$ be Riemannian or pseudo-Riemannian metrics on the same manifold $M^n$. 
Let $L=L(g,\bar g)$ be given by {\rm (\ref{L})}.  

Then,  $g$ and $\bar g$ are geodesically equivalent,  if and only if
\begin{equation} 
\nabla_u L = \frac{1}{2} \bigl(u\otimes l + (u\otimes l)^*\bigr)  \qquad \mbox{for any  tangent vector $u$},
\label{main}
\end{equation}
or, in coordinates,
$$
\nabla_r L^p_q = \frac{1}{2} (\delta^p_r l_q  + g^{ps} l_s g_{rq})
$$
where $\nabla$ is the Levi-Civita connection associated with $g$,  $l = d\,\tr L$   and $C^*$ denotes the operator $g$-adjoint to $C$, i.e., $g(Cu,v)=g(u,C^* v)$.
\end{proposition}

The above  proposition and equation \eqref{main} are due to Sinjukov \cite{sinjukov}, the self-contained proof can also be found in \cite{benenti,EM}.

Since $\bar g$  can be uniquely reconstructed from $g$ and $L$ as 
$\bar g=    
\frac{1}{\det L} g L^{-1}$,  we may replace the pair of metrics $(g, \bar g)$ by the pair $(g, L)$.   For convenience,  we shall say that   a metric $g$ and a $g$-self-adjoint nondegenerate  $(1,1)-$tensor $L$ are  \emph{compatible} if they satisfy (\ref{main}),  so that the compatibility of $g$ and $L$ is just rephrasing  the fact that  $g$ and $\bar g$ are geodesically equivalent.

First of all we recall that the compatibility of $g$ and $L$ implies that the Nijenhuis torsion of $L$ vanishes identically,   see \cite{benenti}.
To make our paper self-contained, we recall some basic facts about the Nijenhuis torsion  and  prove that $N_L \equiv 0$ for any 
$(1,1)$-tensor $L$ satisfying  \eqref{main}.

The \emph{Nijenhuis torsion} of $L$  is the $(1,2)$-tensor field defined by
$$
N_L ( u,v)  = L^2[u,v] - L[Lu,v] - L[u,Lv] + [Lu,Lv]
$$
where $u$ and $v$ are vector fields on $M$ \cite{Haantjes}.
This definition immediately implies
\begin{Lemma}
\label{lem2}
The condition $N_L=0$ admits the following equivalent forms:
\begin{enumerate}
\item  $\mathcal L_{Lu}L - L\mathcal L_uL=0$ for any vector field
$u$,  where $\mathcal L_u$ is the Lie derivative along $u$;

\item $(\nabla_{Lu} L - L\nabla_u L)v$ is symmetric with respect
to $u$ and $v$ for any vector fields $u$ and $v$, where  $\nabla$ is the Levi-Civita connection of an arbitrary  metric $g$ 
(or more generally, any symmetric connection).
\end{enumerate}
\end{Lemma}

\begin{Lemma}[\cite{benenti}]
\label{lem5}
If $L$ satisfies \eqref{main}, then
 $N_L=0$.
\end{Lemma}

{\bf Proof.}  By Lemma \ref{lem2}, we just need to verify that $(\nabla_{Lu} L - L\nabla_u
L)v$ is symmetric with respect to $u$ and $v$. We simply use the
compatibility condition \eqref{main}:
$$
\begin{array}{l}
(\nabla_{Lu} L - L\nabla_u L)v= \frac{1}{2}\bigl( (Lu\otimes l) v
+
(Lu\otimes l)^*v - L (u\otimes l)v  - L(u\otimes l)^*v \bigr) = \\
\\
\frac{1}{2}\bigl( l(v)\cdot Lu + g(Lu,v) \cdot g^{-1}(l)  -
l(v)\cdot Lu  - g(u,v) \cdot L(g^{-1}(l)) \bigr) = \\ \\
\frac{1}{2} \bigl( g(Lu,v) \cdot g^{-1}(l) - g(u,v) \cdot
L(g^{-1}(l)) \bigr)
\end{array}
$$
Here $g^{-1}$ is viewed as the identification map between $T^*_x
M$ and $T_x M$, in particular, $g^{-1}(l)=\mbox{grad}\, \tr L$.
The symmetry of the last expression with respect to $u$ and $v$ is now
evident. \qed 

The next two statements are well known in folklore,  however,  we could not find any reference with
a short proof.

\begin{Lemma}
\label{lem3}
If $N_L=0$, then $N_{p(L)}=0$ for any polynomial $p:\mathbb{R}\to \mathbb{R}$ (with
constant real  coefficients) and, therefore, $N_{f(L)}=0$ for any  function $f:{K}\to \mathbb{C}$ satisfying the assumptions (\ref{I}--\ref{iiiii}) of Section \ref{hijerarchy}.
\end{Lemma}

{\bf Proof.}
The latter statement about  $f(L)$ follows immediately from the definition of $f(L)$.
The proof for a polynomial $p(L)$ is as follows. We use the condition $N_L=0$ in the form
$\mathcal L_{Lu}L = L\mathcal L_uL$ (see Lemma \ref{lem2}).  This identity  implies
$$
\mathcal L_{L^n u} L = \mathcal L_{L(L^{n-1} u)} L = L\mathcal
L_{L^{n-1} u} L =L\mathcal L_{L(L^{n-2} u)} L=L^2\mathcal
L_{L^{n-2} u} L = \dots = L^n\mathcal L_uL,
$$
and, therefore, by linearity
$$
\mathcal L_{p(L)u}L = p(L)\mathcal L_uL.
$$
Thus, we have
$$
(\mathcal L_{p(L)u} - p(L)\mathcal L_u)L=0
$$
Now consider the expression $\mathcal D=\mathcal L_{p(L)u} - p(L)\mathcal
L_u$ as a ``first order differential operator'' which satisfies the obvious 
property $\mathcal D(L^n)= \mathcal D(L^{n-1}) L + L^{n-1} \mathcal D(L)$. Hence, the identity $\mathcal D(L)=0$
immediately implies $\mathcal D(p(L))=0$,~i.e.,
$$
\bigl(\mathcal L_{p(L)u} - p(L)\mathcal L_u\bigr)p(L)=0,
$$
which is exactly the desired condition $N_{p(L)}=0$. \qed

\begin{Lemma}
\label{lem4}
Let $N_L=0$ and $\chi=\chi_1\cdot \chi_2$ be an admissible factorisation of the characteristic polynomial of $L$.
Then there exists a local coordinate system $(x^1,...,x^r,y^{r+1},...,y^n)$ such
that $\partial_{x^1},\dots, \partial_{x^r}$ and
$\partial_{y^{r+1}},\dots, \partial_{y^{n}}$ generate the kernels $D_1$ and $D_2$
of the operators $\chi_1(L)$ and $\chi_2(L)$
respectively and
$$
L(x,y)=\begin{pmatrix} L_1(x) & 0 \\ 0 & L_2(y) \end{pmatrix}.
$$
In particular, the distributions $D_1$ and $D_2$ are integrable.
\end{Lemma}

{\bf Proof.}  At each tangent space $T_xM$
we have the natural decomposition $T_xM=D_1\oplus D_2$  where
$D_i=\ker \chi_i(L)$. 
This decomposition  defines two natural
projectors $P_1$ and $P_2$ onto the subspaces $D_1$ and $D_2$
respectively.  A simple but important observation is that (locally)  these
projectors can be viewed as functions $P_1=f_1(L)$ and
$P_2=f_2(L)$ satisfying the assumptions (\ref{I}--\ref{iiiii}) from Section \ref{hijerarchy}  (see Example \ref{project}).

Thus,  by Lemma \ref{lem3},  $N_{P_i}=0$,  
and we will use this fact to  prove the integrability of 
$D_1$ and $D_2$. 
In terms of the projectors $P_1$ and $P_2$, these distributions can obviously be interpreted as
$D_1=\ker P_2$ and $D_2=\ker P_1$.
 Let $u,v\in D_1=\ker P_2$, then  the condition
$N_{P_2}=0$ gives
$$
P_2^2[u,v]-P_2[P_2u,v]-P_2[u,P_2v]+[P_2u,P_2v]=P^2_2[u,v]=0,
$$
that is, $[u,v]\in\ker P_2^2=\ker P_2=D_1$ which is equivalent to
the integrability of $D_1$ by the Frobenius Theorem. The same is
obviously true for $D_2$ by the same reason.

The integrability of these two distributions is equivalent to the
existence of a coordinate system $(x,y)$ such that
$\partial_{x^1},\dots, \partial_{x^{r}}$ and
$\partial_{y^{r+1}},\dots, \partial_{y^{n}}$ generates the kernels
of the operators $\chi_1(L)$ and $\chi_2(L)$. In
particular, the operator $L$ in this coordinate system has a
block diagonal form:
$$
L(x,y)=\begin{pmatrix} L_1(x,y) & 0 \\ 0 & L_2(x,y) \end{pmatrix}
$$

Whithout loss of generality we may assume that $\det L\ne 0$, otherwise we can (locally)  replace $L$ by $L+ c\cdot {\Id}$  where $c$ is an appropriate constant.

Notice that the operator
$$
L P_1=\begin{pmatrix} L_1(x,y) & 0 \\ 0 & 0
\end{pmatrix},
$$
being a function of $L$,  has zero Nijenhuis torsion.
Thus, for $u=\partial_{y_\alpha}\in D_2=\ker P_1$ we have
$$
\begin{aligned}
\begin{pmatrix} 0 & 0 \\ 0 & 0 \end{pmatrix}=\mathcal L_{LP_1u}(LP_1) &- LP_1\mathcal L_u(LP_1)= LP_1\mathcal
L_u(LP_1)= \\ 
& \begin{pmatrix} L_1 & 0 \\ 0 & 0
\end{pmatrix}
\begin{pmatrix} \partial_{y_\alpha}L_1 & 0 \\ 0 & 0
\end{pmatrix}=
\begin{pmatrix} L_1\partial_{y_\alpha}L_1 & 0 \\ 0 & 0
\end{pmatrix}
\end{aligned}
$$

Since $L_1$ is non-degenerate, we conclude that
$\partial_{y_\alpha}L_1=0$, i.e., $L_1=L_1(x)$. Similarly,
$L_2=L_2(y)$, as needed. \qed

{\bf Proof of Theorems \ref{thm1}, \ref{thm1b}}.  The statements of these theorems are straightforward from Lemmas \ref{lem5},  \ref{lem4}.  
\qed

\subsection{Proof of the generalised Topalov-Sinjukov Theorem \ref{sintop} } \label{topalovsinjukov}

Consider two geodesically equivalent metrics $g$ and $\bar g$ and assume  that $L$ given by \eqref{L} and a function $f:K\to \mathbb{C}$  satisfy the assumptions (\ref{I}--\ref{iiiii}) of Section \ref{hijerarchy}.  Our goal is to prove Theorem \ref{sintop}, i.e., to show that the 
metric $g_f:= g f(L)$ is compatible with  $L$.
To simplify our notation below, we shall denote $g_f$ by $\widetilde g$.

\weg{This statement is an exercise in Tensor Calculus.}

We need to verify that the main equation \eqref{main} for $g$ and $L$ implies the similar relation for
$\widetilde g=g f(L)$ and $L$:
\begin{equation}
 \widetilde \nabla_u L = \frac{1}{2} \left(u\otimes l + (u\otimes
l)^{\widetilde *} \right), \label{eq2}
\end{equation}
where $C^{\widetilde *}$ denotes the operator $\widetilde g$-adjoint of $C$, and $\widetilde \nabla$  is the covariant differentiation with respect to the Levi-Civita connection related to $\widetilde g$.

First of all, we rewrite the condition that we should verify in a
slightly different way by subtracting  \eqref{main} from
(\ref{eq2}):
\begin{equation}
(\widetilde \nabla_u  - \nabla_u) L =\frac{1}{2}  \left((u\otimes
l)^*- (u\otimes l)^{\widetilde *}\right). \label{eq3}
\end{equation}

Now notice that the difference of two covariant derivatives in
the left hand side is not a differential operator, but a tensor
expression of the form
$$ (\widetilde \nabla_u  - \nabla_u) L=T_u L - L
T_u, $$ where $(T_u)^i_j= T_{jk}^i u^k$ and
$T^i_{jk}=\Gamma^i_{jk}-\widetilde \Gamma^i_{jk}$ is a $(1,2)$-tensor
that represents the difference of the two connections related to
$g$  and $\widetilde g$.

Using the obvious fact that $C^{\widetilde *}=f(L)^{-1} C^*
f(L)$ for any operator $C$, we can rewrite (\ref{eq3})  as:
$$
T_u L - L T_u = \frac{1}{2} \left((u\otimes l)^*-
f(L)^{-1}(u\otimes l)^* f(L)\right)
$$

Finally multiplying both sides by $f(L)$ from the left, we see that
the statement of the Topalov--Sinjukov Theorem is equivalent to the
following algebraic relation
\begin{equation}
[f(L) T_u , L]=\frac{1}{2} [f(L),(u\otimes l)^*],
\label{eq5}
\end{equation}
where $[\cdot,\cdot]$ denotes the standard commutator of linear
operators.

To verify it, we compute the tensor $T_u$ explicitly, using the
following equation:
$$
T^s_{jk}\, \widetilde g_{si}+T^s_{ki} \, \widetilde g_{sj}=(\widetilde \nabla_k -
\nabla_k)\, \widetilde g_{ij}=-\nabla_k   \widetilde g_{ij}
$$
This is a system of linear equations w.r.t. $T^s_{jk}\widetilde g_{si}$
with the well-known unique solution:
$$
T^s_{jk}\widetilde g_{si}=\frac{1}{2}  \left( \nabla_j \widetilde g_{sk} +
\nabla_k \widetilde g_{sj} - \nabla_s \widetilde g_{jk}   \right),
$$
or, in invariant terms,
\begin{equation}
\widetilde g(T_u v,w)=\frac{1}{2}  \bigl( (\nabla_u \widetilde g)(v,w) +
(\nabla_v \widetilde g) (u,w)  - (\nabla_w \widetilde g) (u,v) \bigr)
\label{eq6}
\end{equation}

In our particular case,  we have:
\begin{equation}
\begin{array}{l}
\nabla_u \widetilde g = g \nabla_u f(L) = g f'_{t=0}(L + t\nabla_u L)
=  \\
\\   f'_{t=0} \left(L + t \frac{1}{2}\left(u\otimes l +
(u\otimes
l)^*\right)\right)=\\
\\
\frac{1}{2} g  f'_{t=0}(L + t (u\otimes l))  + \frac{1}{2} g
f'_{t=0}(L + t(u\otimes l)^*)=\frac{1}{2}g (A_u + A_u^*),
\end{array}
\label{eq7}
\end{equation}
where $A_u=f'_{t=0}(L + t (u\otimes l))$ and $A_u^*=f'_{t=0}(L + t
(u\otimes l)^*)$. The operators $A_u$ and $A_u^*$ are $g$-adjoint to each other
and satisfy the following important property.

\begin{Lemma} In the notation above, 
$A^*_u v = A^*_v u$.
\end{Lemma}

{\bf Proof.} If $f(L)=L$ then $A^*_u=(u\otimes
l)^*=g^{-1}(l)\otimes g(u)$ and we have $A^*_u v=g(u,v)
g^{-1}(l)=A^*_v u$ (here we consider $g$ and $g^{-1}$ as
identification operators between the tangent and cotangent
spaces).

If $f(L)=L^2$, then $A^*_u =L(u\otimes l)^* + (u\otimes l)^* L$
and $A^*_u v=g(u,v) \cdot Lg^{-1}(l) + g(u,Lv) \cdot
g^{-1}(l)=A^*_v u$ (we use, of course, the fact that $L$ is
$g$-self-adjoint).

More generally, for $f(L)=L^n$, we have $$A^*_u v=\sum_{m=0}^{n-1}
g(L^{n-m}u,v)\cdot L^{m-1} g^{-1}(l)=A^*_v u.$$

Thus, the statement holds for any polynomial $f(L)$ and,
therefore, for any smooth function $f$ satisfying (\ref{I}--\ref{iiiii}) in Section \ref{hijerarchy}.

Using these properties, we get a surprisingly simple result by
substituting (\ref{eq7}) into (\ref{eq6}):
$$
\begin{array}{l}
\widetilde g(T_u v,w)=\frac{1}{4} \bigl( g((A_u\!+\!A^*_u) v,w)+
g((A_v\!+\!A^*_v) u,w) \\
\\- g((A_w\!+\!A^*_w) u,v)\bigr)=\frac{1}{2}g(A^*_u v,w)
\end{array}
$$
In other words, $f(L) T_u = A_u^* = \frac{1}{2} f'_{t=0}(L + t
(u\otimes l)^*_g)$.

We are now ready to complete the proof. We use the following
simple matrix relation which obviously holds for any $L,B$ and
$f$:
$$
0=[f(L+tB), L+tB]'_{t=0}=[f'_{t=0}(L+tB), L] + [f(L),B],
$$
that is,
$$
[f'_{t=0}(L+tB), L] = [B,f(L)].
$$
In our case $B=(u\otimes l)^*_g$ and $f'_{t=0}(L+tB)=2f(L)T_u$
which gives exactly (\ref{eq5}), as required. \qed

\subsection{Proof of Theorem \ref{thm2}}

Theorems \ref{thm1}, \ref{thm1b}   reduce our consideration to the block-diagonal case in the sense that
 the metric $g$ and  operator $L$  both have block-diagonal form in local coordinates
$(x, y)=(x^1,...,x^r,y^{r+1},...,y^n)$, more precisely
\begin{equation}
g=\begin{pmatrix}
g_1(x,y) & 0 \\  0 & g_2(x,y)
\end{pmatrix}, \qquad
L=\begin{pmatrix}
L_1(x) & 0 \\  0 & L_2(y)
\end{pmatrix}.
\label{blockdiag}
\end{equation}

 Notice first that in such a situation, the main equation \eqref{main}  has a rather special form.  Namely, it can naturally be divided into three parts each of which has its own meaning and can be treated separately (up to some extent).

{\bf Index notation convention}.
For convenience, throughout the rest of the  paper we denote the indices for $x^1, \dots, x^r$ by
Latin letters $i,j,k,l,m=1, \dots, r$,  and those for $y^{r+1}, \dots  y^n$ by Greek letters $\alpha, \beta, \gamma, \delta = r+1, \dots, n$.  The indices $p,q,r,s,t$ will serve for both cases, i.e., 
$p,q,r,s,t=1,\dots, n$.

Let $u \in D_1$, then  the main equation \eqref{main}  can be written in a block form as
$$
 \begin{pmatrix}
\nabla_u L^i_j  & \nabla_u L^i_\alpha \\
\nabla_u L^\beta_j & \nabla_u L^\beta_\alpha
\end{pmatrix} =
\begin{pmatrix}
\frac{1}{2} \left(u^i
\frac{\partial \tr L}{\partial x^j} +  g^{il} \frac{\partial \tr L}{\partial x^l}u^m  g_{mj} \right) &\frac{1}{2} u^i
\frac{\partial \tr L}{\partial y^\alpha}\\  \frac{1}{2}g^{\beta\gamma} \frac{\partial \tr L}{\partial y^\gamma} u^m g_{mj}   & 0
\end{pmatrix}
$$

We rewrite it for each block separately taking into account the block-diagonal form of $g$ and $L$ and the fact that $\frac{\partial \tr L}{\partial x^l}=\frac{\partial \tr L_1}{\partial x^l}$ and $\frac{\partial \tr L}{\partial y^\alpha}=\frac{\partial \tr L_2}{\partial y^\alpha}$:
 \begin{equation}
 u^k \frac{\partial (L_1)^i_j}{\partial x^k} + u^ k \Gamma^{i}_{km}
(L_1)^m_j -u^k \Gamma^m_{kj} (L_1)^i_m =\frac{1}{2} \left(u^i
\frac{\partial \tr L_1}{\partial x^j} +  (g_1)^{il} \frac{\partial \tr L_1}{\partial x^l}u^m  (g_1)_{mj} \right)
\label{main1}
 \end{equation}

\begin{equation}
 u^k\Gamma^j_{k\beta} (L_2)^\beta_\alpha
-u^k\Gamma^m_{k\alpha} (L_1)^j_m =\frac{1}{2} u^j  \frac{\partial \tr L_2}{\partial y^\alpha}
 \label{main2}
 \end{equation}

\begin{equation}
 u^k \Gamma_{\gamma k}^\alpha (L_2)^\gamma_\beta -
(L_2)^\alpha_\gamma u^k \Gamma^\gamma_{\beta k} = 0
\label{main3}
 \end{equation}

We omit the equation for  $\nabla_u L^\beta_j$ because it is obtained from \eqref{main2} by using the fact that the left hand side and right hand side on \eqref{main} are both $g$-self-adjoint).

First of all, we notice that the Christoffel  symbols $\Gamma^i_{jk}$  ($i,j,k =1,\dots , r<n$) of the Levi-Civita connection associated with $g$ coincide with those for the metric $g_1$ defined on the leaves of the integrable distribution $D_1$.  Thus,  the first equation \eqref{main1}  simply means that  $g_1(x,y)$ and $L_1(x)$ are compatible on each leaf of $D_1$, i.e., for every fixed $y$. 

To simplify \eqref{main2} and \eqref{main3}, we notice that 
$$
\Gamma_{\alpha k}^i = \frac{1}{2} g^{im}\frac{\partial g_{mk}}{\partial y^\alpha}  \qquad \mbox{and} \qquad
\Gamma_{\alpha k}^\beta=\frac{1}{2} g^{\beta\gamma} \frac{\partial g_{\gamma\alpha}}{\partial x^k}
$$

Substituting these expressions into  \eqref{main2} and \eqref{main3}  and taking into account that $u\in  D_1$  is an arbitrary vector,  we obtain
\begin{equation}
g^{jm}\frac{\partial g_{mk}}{\partial y^\beta} 
(L_2)^\beta_\alpha
- (L_1)^j_m g^{ml}\frac{\partial g_{lk}}{\partial y^\alpha} 
 =  \delta^j_k  \frac{\partial \tr L_2}{\partial y^\alpha}
\label{main22}
\end{equation}
and 
\begin{equation}
 g^{\alpha\beta} \frac{\partial g_{\beta\gamma}}{\partial x^k}
   (L_2)^\gamma_\beta -
(L_2)^\alpha_\gamma  
g^{\gamma\beta} \frac{\partial g_{\beta\alpha}}{\partial x^k}
 = 0
\label{main32}
 \end{equation}

Summarizing this discussion and rewriting \eqref{main22} and \eqref{main32} in a shorter ``matrix'' form, we obtain

\begin{proposition}
\label{prop1}
For $g$ and $L$ of block-diagonal form  \eqref{blockdiag},  the
compatibility equation \eqref{main} for $u\in D_1$ is equivalent
to the following three conditions:

{\rm Condition 1: }  $g_1(x,y)$ and $L_1(x)$ are compatible on each leaf of $B_1$  (here $x$ are coordinates on  a leaf, $y$ is a parameter which determines  the leaf);

{\rm Condition 2:}
\begin{equation}
 (g_1^{-1} d_y g_1) \,  L_2 - L_1 \, (g_1^{-1} d_y g_1)  =  \Id_{D_1}
\otimes d_y \tr L_2;
\label{cond2}
\end{equation}

{\rm Condition 3:}
\begin{equation}
(g_2^{-1} d_x g_2) \, L_2 - L_2 \, (g_2^{-1} d_x g_2)  =0.
\label{cond3}
\end{equation}

\end{proposition}

As we shall see now, Theorems \ref{thm2}  and \ref{thm3}  are
algebraic corollaries of these equations, generalized  Topalov-Sinjukov
Theorem \ref{sintop} and vanishing of $N_L$.

The first part of  Theorem \ref{thm2}  states that  the metric
$$
h=\begin{pmatrix}
h_1 & 0 \\ 0 & h_2
\end{pmatrix}=\begin{pmatrix}  g_1 \chi_2(L_1)^{-1} & 0 \\ 0 & g_2 \chi_1(L_2)^{-1} \end{pmatrix}
$$
is of local product structure. In other words, we need to check that
$$
\dfrac{\partial (h_1)_{lm}}{\partial y^\alpha}=0
\quad\mbox{and}\quad \dfrac{\partial
(h_2)_{\alpha\beta}}{\partial x^m}=0.
$$
 We shall verify the first condition only, the proof for the  second is similar. We start with straightforward computation.

\begin{Lemma}
Let $h_1 =  g_1\, \chi_2 (L_1)^{-1} $, then
\begin{equation}
\dfrac{\partial (h_1)_{lm}}{\partial y_\alpha}= g_{li} \left(
g^{ij}\dfrac{ \partial g_{jk}}{\partial y^\alpha}   - \bigl(\chi_2(L_1)^{-1}\bigr)^i_j \dfrac{\partial
\chi_2(L_1)^j_k}{\partial y^\alpha}  \right)
\bigl(\chi_2(L_1)^{-1}\bigr)^k_m 
\label{h11}
\end{equation}
or, in coordinate-free form, 
\begin{equation}
d_y  h_1 = g_{1}  \left( g^{-1}_1  d_y g_1 -
\chi_2(L_1)^{-1}  d_y {\chi_2(L_1)}  \right) \chi_2(L_1)^{-1}.
\label{h112}
\end{equation}
\label{h1}
\end{Lemma}

{\bf Proof.} Thinking of $h_1=g_1\, \chi_2 (L_1)^{-1} $ as a
product of two $r\times r$-matrices,  we have
$$
\begin{array}{c}
\dfrac{\partial}{\partial y_\alpha}  h_1 =
\dfrac{\partial}{\partial y^\alpha} \bigl( g_1 \chi_2(L_1)^{-1}\bigr)=
\dfrac{\partial g_1}{\partial y^\alpha}  \, \chi_2(L_1)^{-1} +
g_1 \dfrac{\partial \chi_2(L_1)^{-1}}{\partial y^\alpha}  = \\   \\
\dfrac{\partial g_1}{\partial y^\alpha}  \, \chi_2(L_1)^{-1} -
g_1 \, \chi_2(L_1)^{-1} \dfrac{\partial \chi_2 (L_1)}{\partial y^\alpha}  \chi_2(L_1)^{-1} = \\ \\
g_{1}  \left( g^{-1}_1 \dfrac{ \partial g_1}{\partial y^\alpha} -
\chi_2(L_1)^{-1} \dfrac{\partial \chi_2(L_1)}{\partial
y^\alpha}  \right) \chi_2(L_1)^{-1},
\end{array}
$$
as needed. \qed 

Thus, we need to prove that the expression in brackets in the left hand side of \eqref{h11} (or   \eqref{h112} ) is actually zero.
To do so, we need some properties of the differential of the characteristic polynomial of $L$ in the case when  $N_L$ vanishes.

\begin{Lemma} \label{lemm} 
Let $L$ be a tensor with zero Nijenhuis torsion,  $l = d\,
\tr L$  be the differential of $\tr L$ viewed as a covector and $\chi(t)=\det (t\cdot \Id-L)$ the characteristic polynomial of $L$ viewed as a smooth function on $M$ with $t$ as a formal parameter.  Then the differential of $\chi(t)$  satisfies the following relation:
\begin{equation}
d \chi(t)\, L - t \cdot d \chi(t)=\chi(t) \cdot l
\label{eq10}
\end{equation}
\end{Lemma}

\begin{Rem}  {\rm  In  (\ref{eq10}),    the right and left hand sides are both covectors, i.e. elements of the cotangent space.  The expression  $d \chi(t)\, L$ means that we apply the operator $L$ to the covector  $d\chi$ using right multiplication.  In coordinates, this means 
$\bigl(d \chi(t)\, L\bigr)_s= \dfrac{\partial \chi(t)}{\partial x^p} L^p_s$. The multiplication denoted by $\cdot$  simply means  multiplying a covector by a scalar function.}
\end{Rem}

{\bf Proof of Lemma \ref{lemm}.}  The differential of the polynomial
$$
\chi(t)=\det (t\cdot \Id-L)= a_0  + a_1 t + a_2 t^2 + \ldots + t^n,
$$
is defined to be
$$
d\chi(t)= d a_0 + t \cdot da_1  + t^2 \cdot da_2 + \ldots +  t^{n-1} \cdot da_{n-1}.
$$
Instead of differentiating each coefficient separately, we are
going to differentiate the whole polynomial  (thinking of $t$ as
a certain constant).  For the operators $L$ with $N_L=0$,   the following
property is well known (see, for example, \cite[Lemma1]{benenti}):
$$
\bigl(  d \log|\det  L| \bigr)  \,  L  = d\, \tr \, L= l
$$ 
Notice that  $N_{t\cdot \Id-L}=0$ and apply the above identity to  $\chi(t)=\det (t\cdot \Id-L )$:
$$
\bigl(  d \log |\chi(t)| \bigr)  \, (t\cdot \Id-L ) = d\, \tr (t\cdot \Id-L )=-d\, \tr L =  -l.
$$
Thus
$
 \frac{1}{\chi(t)} \, d  \chi(t)   \, (L-t\cdot \Id)  =  l,
$
which is equivalent to (\ref{eq10}).  \qed

In fact, we need to compute the differential of the characteristic
polynomial in the case when $L=L_2$ and $t$ is replaced by $L_1$.
More precisely we shall need the expression:
$$
\bigl(d_y \chi_2 (L_1)\bigr)^i_{j,\alpha}=\frac{\partial}{\partial y^\alpha} (\chi_2
(L_1))^i_j
$$

This object can be naturally viewed as
an element of the space $\mbox{Hom}(D_1, D_1)\otimes D_2^*$.
Since (\ref{eq10}) is purely algebraic (in the sense that the
nature of $t$ is not important) and $t=L_1$ does not depend on  $y$,  we can reformulate (\ref{eq10}) for our  special case as follows:
\begin{equation}
d_y \chi_{2} (L_1) \, L_2 - L_1 \, d_y \chi_{2} (L_1) = \chi_{2} (L_1) \otimes l_2
\label{eq12}
\end{equation}
or, in coordinates,
$$
\bigl(d_y \chi_{2} (L_1)\bigr)^i_{j,\beta} (L_2)^\beta_\alpha - (L_1)^i_k \bigl( d_y \chi_{2} (L_1)\bigr)^k_{j,\alpha} =
(\chi_{2} (L_1))^i_j \, (l_2)_\alpha,
$$
where  $l_2=d \,\tr L_2$.

Multiplying the both sides of \eqref{eq12} by $\chi_2 (L_1)^{-1}$, and using that $L_1$ commute with $\chi_2(L_1)$,   we obtain

\begin{Lemma} If $N_L=0$, then the expression $\chi_2 (L_1)^{-1} \, d_y \chi_{2} (L_1)$ satisfies the following relation  
\begin{equation}
\big( \chi_2 (L_1)^{-1} \, d_y \chi_{2} (L_1)\big)  \, L_2 - L_1 \, \big( \chi_{2} (L_1)^{-1} \,  d_y \chi_{2} (L_1) \big)= {\Id}_{D_1}
\otimes  l_2
\label{eq14}
\end{equation}
or, in coordinates,
$$
(\chi_{2} (L_1)^{-1})^j_m \bigl(d_y \chi_{2} (L_1)\bigr)^m_{k,\beta} (L_2)^\beta_\alpha -
(L_1)^j_m (\chi_{2} (L_1)^{-1})^m_l  \bigl(d_y \chi_{2} (L_1)\bigr)^l_{k,\alpha} = \delta^j_k \,
 (l_2)_\alpha.
$$
\label{lemma7}
\end{Lemma}

We are now ready to complete the proof of Theorem \ref{thm2}.
Let $L$ and $g$ be compatible.  Then Condition 2 from Proposition \ref{prop1} holds.  Comparing this condition \eqref{cond2}  with the relation \eqref{eq14},   we see that  $\chi_2 (L_1)^{-1} \, d_y \chi_{2} (L_1)$ and $g^{-1}_1  d_y g_1$ satisfy the same relation.
If we fix index $k$  in  \eqref{cond2}  and  \eqref{eq14}, we shall see that these relations can be treated as a linear matrix equation of the form
$$
X \, L_2  - L_1 \, X  = C  
$$
Where $X$ is an (unknown) matrix of dimension $r\times (n-r)$,  and $L_1$, $L_2$, $C$ are given matrices of dimensions $r\times r$, $(n-r)\times (n-r)$ and $r\times(n-r)$ respectively.

 It is a simple fact in Linear Algebra that if $L_1$ and
$L_2$ have no common eigenvalues, then the solution to this equation
is unique for any $C$. Hence we conclude that 
$$\chi_2 (L_1)^{-1} \, d_y \chi_{2} (L_1) = g^{-1}_1  d_y g_1,$$
which  immediately implies 
$\frac{\partial (h_1)_{ij}}{\partial y^\alpha} = 0$ (see Lemma \ref{h1}).

Analogously  (i.e., just by interchanging the distributions $D_1$ and $D_2$),  we can check that $\frac{\partial (h_2)_{\alpha\beta}}{\partial x^i} = 0$. This means that  in the metric
$$
h = \begin{pmatrix} h_1 & 0 \\ 0  & h_2 \end{pmatrix}
$$ the block $h_1$ depends on the $x-$coordinates only, and  the block $h_2$ depends on the $y-$coordinates only, i.e.,  $(h, B_1, B_2)$ is  
a  local product structure.  Thus, the first  statement of  Theorem \ref{thm2} is proved.

The second statement of the theorem \ref{thm2}  says that  the restrictions of $h$ and $\bar h$ on the same leaf of $B_1$, i.e., the metrics $h_1(x)=g_1 \chi_2(L_1)^{-1}$ and $\bar h_1(x)= 
\tfrac{1}{\chi_2(0)} g_1 \chi_2(L_1)^{-1}$   are geodesically equivalent.  Equivalently, we can reformulate this saying that $h_1$ and $L_1$ are compatible on each leaf of the foliation $F_1$.  To prove this fact we only need to notice that  on each fixed leaf
$\chi_2(L_1)$ is a polynomial (with constant coefficients) in $L_1$ and, therefore,  $\chi_2(L_1)^{-1}$ is a ``good'' function of $L_1$ so that the compatibility of
$h_1= g_1 \chi_2(L_1)^{-1}$ and $L_1$  follows from the Theorem  \ref{sintop} and Condition 1 of   Proposition \ref{prop1}.

\subsection{Proof of Theorem \ref{thm3} }

Now we are going to show that the compatibility of  the pairs $h_1(x),L_1(x)$ and $h_2(y), L_2(y)$ imply  the compatibility of 
$$
g=\begin{pmatrix} g_1(x,y) & 0 \\
0 & g_2(x,y)
   \end{pmatrix}=\begin{pmatrix}
   h_1 \chi_2(L_1) & 0 \\ 0 & h_2 \chi_1(L_2)
   \end{pmatrix}
 \qquad \mbox{and} \qquad  L = \begin{pmatrix}  L_1 &  0 \\ 0 & L_2 \end{pmatrix},
$$
which is equivalent to the statement of Theorem \ref{thm3}.

We shall verify the compatibility condition \eqref{main}  for $u\in D_1$  only  (the proof for $u\in D_2$ is absolutely similar).
Since $g$ and $L$ are of block-diagonal form \eqref{blockdiag}, we may use Proposition \ref{prop1} and replace  \eqref{main} by Conditions 1--3.

Condition 1 (i.e., compatibility of $g_1=h_1 \chi_2(L_1)$  and $L_1$ on leaves of $B_1$) immediately follows from the Topalov-Sinjukov theorem and compatibility of $h_1$ and $L_1$.

To verify Condition 2,  we notice that  Lemma \ref{h1} and the condition $\frac{\partial}{\partial y^\alpha} h_1(x) = 0$  imply the relation:
$$
g^{-1}_1  d_y g_1 -
\chi_2(L_1)^{-1}  d_y {\chi_2(L_1)}=0 
$$

By  Lemma \ref{lemma7}, $ \chi_2(L_1)^{-1}  d_y {\chi_2(L_1)} $  satisfies \eqref{eq14},  therefore so does  $g^{-1}_1  d_y g_1$, which gives exactly Condition 2.

The last  relation \eqref{cond3}  (Condition 3) is the matrix equation
$$
\left[L_2,  g_2^{-1} \frac{\partial g_2}{\partial x^k}\right]=0,  \quad\mbox{for any $k=1,\dots,r$.}
$$

By our construction,   $g_2= h_2 \chi_1(L_2)$, where
$\chi_1(L_2)$ is a polynomial of the form $\chi_1(L_2)=
a_0(x) + a_1(x) L_2 + a_2(x)L^2_2+ a_3(x) L^3_2 + \dots$ and
neither $h_2$ nor $L_2$ depend on $x^k$. Hence
$$
\begin{array}{l}
g_2^{-1}\dfrac{\partial g_2}{\partial x^k} = \chi_1(L_2)^{-1}
h_2^{-1} h_2 \left( \dfrac{\partial a_0}{\partial x^k} +
\dfrac{\partial a_1}{\partial x^k}  L_2 + \dfrac{\partial
a_2}{\partial x^k} L^2_2+ \dfrac{\partial a_3}{\partial x^k} L^3_2
+ \dots \right)=\\
\chi_1(L_2)^{-1}\left( \dfrac{\partial a_0}{\partial x^k} +
\dfrac{\partial a_1}{\partial x^k}  L_2 + \dfrac{\partial
a_2}{\partial x^k} L^2_2+ \dfrac{\partial a_3}{\partial x^k} L^3_2
+ \dots \right)=f(L_2)
\end{array}
$$
where $f(L_2)$ is a  function of $L_2$ (depending on $x$ as
a parameter). Thus, (\ref{cond3}) holds since $[L_2, f(L_2)]=0$ for  
every function $f$ (satisfying assumptions (\ref{I}--\ref{iiiii}) of Section \ref{hijerarchy}).
This completes the proof of Theorem \ref{thm3}. \qed

{\em Acknowledgement:} The first author thanks I~Zakharevich for useful discussions.
The second author thanks 
Deutsche Forschungsgemeinschaft (Priority Program 1154 --- Global Differential Geometry) for partial financial support, and V. Shevchishin and  N. Higham for useful discussions.


\end{document}